\documentclass[11pt]{amsart}
\usepackage{amssymb}
\usepackage{amsmath}
\usepackage{graphicx,psfrag}

\title{A survey on Hamilton cycles in directed graphs}

\date{\today}
\author{Daniela K\"uhn and Deryk Osthus}

\newtheorem{firstthm}{Proposition}
\newtheorem{theorem}[firstthm]{Theorem}

\newtheorem{lemma}[firstthm]{Lemma}

\newtheorem{conj}[firstthm]{Conjecture}
\newtheorem{question}[firstthm]{Question}

\def\noproof{{\unskip\nobreak\hfill\penalty50\hskip2em\hbox{}\nobreak\hfill%
       $\square$\parfillskip=0pt\finalhyphendemerits=0\par}\goodbreak}
\def\endproof{\noproof\bigskip}
\newdimen\margin   
\def\textno#1&#2\par{%
   \margin=\hsize
   \advance\margin by -4\parindent
          \setbox1=\hbox{\sl#1}%
   \ifdim\wd1 < \margin
      $$\box1\eqno#2$$%
   \else
      \bigbreak
      \hbox to \hsize{\indent$\vcenter{\advance\hsize by -3\parindent
      \it\noindent#1}\hfil#2$}%
      \bigbreak
   \fi}

\begin{document}
\maketitle

\def\COMMENT#1{}
\def\TASK#1{}

\def\eps{{\varepsilon}}
\newcommand{\ex}{\mathbb{E}}
\newcommand{\pr}{\mathbb{P}}
\newcommand{\C}{\mathcal{C}}
\newcommand{\M}{\mathcal{M}}
\begin{abstract} \noindent
We survey some recent results on long-standing conjectures regarding 
Hamilton cycles in directed graphs, oriented graphs and tournaments.
We also combine some of these to prove the following approximate result
towards Kelly's conjecture on Hamilton decompositions of regular tournaments: the edges of every regular tournament
can be covered by a set of Hamilton cycles which are `almost' edge-disjoint.
We also highlight the role that the notion of `robust expansion' plays in several of the proofs.
New and old open problems are discussed.
\end{abstract}

\section{Introduction} \label{intro}

The decision problem of whether 
a graph has a Hamilton cycle is NP-complete and so a satisfactory characterization of Hamiltonian
graphs seems unlikely. Thus it makes sense to ask for degree
conditions which ensure that a graph has a Hamilton cycle.
One such result is Dirac's theorem~\cite{dirac}, which states that every graph on $n \ge 3$ vertices with 
minimum degree at least $n/2$ contains a Hamilton cycle. This is strengthened by 
Ore's theorem~\cite{ore}:
If $G$ is a graph with $n \ge 3$ vertices such that every pair $x\neq y$
of non-adjacent vertices satisfies $d(x)+d(y) \ge n$, then $G$ has a Hamilton cycle.
Dirac's theorem can also be strengthened considerably by allowing many of the vertices to have small degree:
P\'osa's theorem~\cite{posa} states that a graph on~$n\ge 3$
vertices has a Hamilton cycle if its degree sequence $d_1\le  \dots \le d_n$ satisfies $d_i \geq i+1$ 
for all $i<(n-1)/2$ and if additionally $d_{\lceil n/2\rceil} \geq \lceil n/2\rceil$ when~$n$ is odd. 
Again, this is best possible -- none of the degree conditions can be relaxed.
Chv\'atal's
theorem~\cite{chvatal} is a further generalization. It characterizes all those degree sequences which
ensure the existence of a Hamilton cycle in a graph: suppose that the degrees of the graph $G$
are $d_1\le \dots \le d_n$. If $n \geq 3$ and $d_i \geq i+1$ or $d_{n-i} \geq n-i$
for all $i <n/2$ then $G$ is Hamiltonian. This condition on the degree sequence is
best possible in the sense that for any degree sequence $d_1\le \dots \le d_n$ violating this condition there
is a corresponding graph with no Hamilton cycle whose degree sequence dominates $d_1,\dots ,d_n$.
These four results are among the most general and well-known Hamiltonicity conditions.
There are many more -- often involving additional structural conditions like planarity.
The survey~\cite{gould} gives an extensive overview (which concentrates on undirected graphs).

In this survey, we concentrate on recent progress for directed graphs.
Though the problems are equally natural for directed graphs, it is usually much more difficult to obtain satisfactory results.
Additional results beyond those discussed here can be found in the corresponding chapter of the monograph~\cite{digraphsbook}.
In Section~\ref{digraph}, we discuss digraph analogues and generalizations of the above four results.
The next section is devoted to oriented graphs -- these are obtained from undirected graphs by orienting the
edges (and thus are digraphs without $2$-cycles).
Section~\ref{tournament} is concerned with tournaments.
Section~\ref{general} is devoted to several generalizations of the notion of a Hamilton cycle,
e.g.~pancyclicity and  $k$-ordered Hamilton cycles.
The final section is devoted to the concept of `robust expansion'.
This has been useful in proving many of the recent results discussed in this survey.
We will give a brief sketch of how it can be used.
In this paper, we also use this notion (and several results from this survey)
to obtain a new result (Theorem~\ref{cover}) which gives further support to Kelly's conjecture on
Hamilton decompositions of regular tournaments. In a similar vein, we use a result of~\cite{cko}
to deduce that the edges of every sufficiently dense regular (undirected) graph can be covered
by Hamilton cycles which are almost edge-disjoint (Theorem~\ref{RegularII}).

\section{Hamilton cycles in directed graphs} \label{digraph}

\subsection{Minimum degree conditions}

For an analogue of Dirac's theorem in directed graphs it is
natural to consider the \emph{minimum semidegree~$\delta^0(G)$} of a digraph $G$,
which is the minimum of its minimum outdegree~$\delta^+(G)$ and its minimum
indegree~$\delta^-(G)$. (Here a directed graph may have two edges between a pair of vertices, 
but in this case their directions must be opposite.)
The corresponding result is a theorem of Ghouila-Houri~\cite{gh}.
\begin{theorem}[Ghouila-Houri~\cite{gh}] \label{hamGH}
Every strongly connected digraph on $n$ vertices with $\delta^+(G)+\delta^-(G) \ge n$ contains a Hamilton cycle.
In particular, every digraph with $\delta^0(G)\ge n/2$ contains a Hamilton cycle. 
\end{theorem}
(When referring to paths and cycles in directed
graphs we usually mean that these are directed, without mentioning this explicitly.)

For undirected regular graphs, Jackson~\cite{jacksonreg} 
showed that one can reduce the degree condition
in Dirac's theorem considerably if we also impose a connectivity condition, i.e.~every $2$-connected $d$-regular graph on $n$
vertices with $d\ge n/3$ contains a Hamilton cycle. Hilbig~\cite{hilbig} improved the degree condition
to~$n/3-1$ unless $G$ is the Petersen graph or another exceptional graph.
The example in Figure~1 shows that the degree condition cannot be reduced any further.
Clearly, the connectivity condition is necessary.
\begin{figure}\label{digrapregular}
\centering\footnotesize
\includegraphics[scale=0.4]{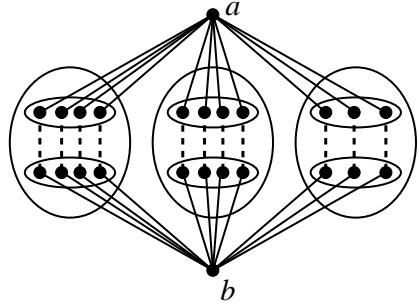}
\caption{A $(3s-1)$-regular $2$-connected graph $G$ on $n=9s+2$ vertices with no Hamilton cycle.
To construct $G$, start with $3$ disjoint cliques on $3s$ vertices each.
In the $i$th clique choose disjoint sets $A_i$ and $B_i$ with $|A_i|=|B_i|$ and
$|A_1|=|A_2|=s$ and $|A_3|=s-1$. Remove a perfect matching between $A_i$ and $B_i$ for 
each $i$. Add 2 new vertices $a$ and $b$, where $a$ is connected to all vertices in the sets $A_i$
and $b$ is connected to all vertices in all the sets $B_i$.}
\end{figure}
We believe that a similar result should hold for directed graphs too.
\begin{conj} \label{diregconj}
Every strongly $2$-connected $d$-regular digraph on $n$ vertices with $d\ge n/3$ contains a
Hamilton cycle.
\end{conj}
Replacing each edge in Figure~1 with two oppositely oriented edges shows that the
degree condition cannot be reduced. Moreover, it is not hard to see that the strong
$2$-connectivity cannot be replaced by just strong connectivity.%
   \COMMENT{Start with 2 complete digraphs $A$ and $B$ on $n/2$ vertices.
Delete the edges of an $(n/4-1)$-cycle $C_A$ in $A$ and of
an $n/4$-cycle $C_B$ in $B$. Add a new vertex $x$ and join it to all
vertices on $C_A$ and $C_B$ with edges in both directions. Then the resulting
digraph is $(n/2-1)$-regular and strongly connected but does not contain
a Hamilton cycle.} 

\subsection{Ore-type conditions}

Woodall proved the following digraph version of Ore's theorem, which generalizes Ghouila-Houri's theorem.
$d^+(x)$ denotes the outdegree of a vertex $x$, and $d^-(x)$ its indegree.
\begin{theorem}[Woodall~\cite{woodall}] \label{woodall}
Let $G$ be a strongly connected digraph on $n \ge 2$ vertices.
If $d^+(x)+d^-(y)\ge n$ for every pair $x\neq y$ of vertices for which there is no edge from $x$ to $y$,
then $G$ has a Hamilton cycle.
\end{theorem}
Woodall's theorem in turn is generalized by Meyniel's theorem, where the degree condition is
formulated in terms of the total degree of a vertex. Here the \emph{total degree $d(x)$ of $x$} is defined as
$d(x):=d^+(x)+d^-(x)$.
\begin{theorem}[Meyniel~\cite{meyniel}]
Let $G$ be a strongly connected digraph on $n \ge 2$ vertices.
If $d(x)+d(y) \ge 2n-1$ for all pairs of non-adjacent vertices
in $G$, then $G$ has a Hamilton cycle.
\end{theorem}
The following conjecture of Bang-Jensen, Gutin and Li~\cite{BGL} would strengthen Meyniel's theorem by 
requiring the degree condition only for dominated pairs of vertices%
    \COMMENT{can replace dominated by dominating} 
(a pair of vertices is \emph{dominated} if there is a vertex which sends an edge to both of them).
\begin{conj}[Bang-Jensen, Gutin and Li~\cite{BGL}] \label{BGLconj}
Let $G$ be a strongly connected digraph on $n \ge 2$ vertices.
If $d(x)+d(y) \ge 2n-1$ for all dominated pairs of non-adjacent vertices
in $G$, then $G$ has a Hamilton cycle.
\end{conj}
\begin{figure}\label{digraphF}
\centering\footnotesize
\includegraphics[scale=0.45]{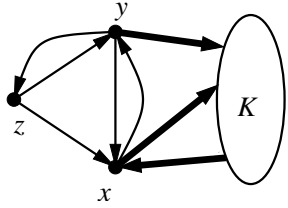}
\caption{An extremal example for Conjecture~\ref{BGLconj}:
let $F$ be the digraph obtained from the complete digraph $K=K_{n-3}^{\leftrightarrow}$
and a complete digraph on 3 other vertices $x,y,z$ as follows:
remove the edge from $x$ to $z$, add all edges in both directions between $x$ and $K$
and all edges from $y$ to $K$.}
\end{figure}
An extremal example $F$ can be constructed as in Figure~2.
To see that $F$ has no Hamilton cycle, note that every Hamilton path in $F-x$ has to start at $z$.
Also, note that the only non-adjacent (dominated) pairs of vertices are $z$ together with a vertex $u$ in $K$
and these satisfy $d(z)+d(u)=2n-2$.

Some support for the conjecture is given e.g.~by the following result of Bang-Jensen, Guo and Yeo~\cite{BGY}:
if we also assume the degree condition for all pairs of non-adjacent vertices which have a common outneighbour, then
$G$ has a $1$-factor, i.e.~a union of vertex-disjoint cycles covering all the vertices of $G$.

There are also a number of degree conditions which
involve triples or 4-sets of vertices, see e.g.~the corresponding chapter in~\cite{digraphsbook}.

\subsection{Degree sequences forcing Hamilton cycles in directed graphs}

Nash-Williams~\cite{nw}
raised the question of a digraph analogue of Chv\'atal's theorem quite soon after the
latter was proved:
for a digraph~$G$ it is natural to consider both its outdegree sequence $d^+ _1,\dots , d^+ _n$
and its indegree sequence $d^- _1,\dots , d^- _n$. Throughout, we take the convention that
$d^+ _1\le \dots \le d^+ _n$ and $d^- _1 \le \dots \le  d^- _n$ without mentioning this explicitly.
Note that the terms $d^+ _i$ and $d^- _i$ do not necessarily correspond to the degree of the same vertex
of~$G$. 
\begin{conj}[Nash-Williams~\cite{nw}]\label{nw}
Suppose that $G$ is a strongly connected digraph on $n \geq 3$ vertices
such that for all $i < n/2$
\begin{itemize}
\item[{\rm (i)}]  $d^+ _i \geq i+1 $ or $ d^- _{n-i} \geq n-i $,
\item[{\rm (ii)}] $ d^- _i \geq i+1$ or $ d^+ _{n-i} \geq n-i.$
\end{itemize}
Then $G$ contains a Hamilton cycle.
\end{conj}
It is even an open problem whether the conditions imply the existence of a cycle through
any pair of given vertices (see~\cite{bt}).  
The following example  shows that the 
degree condition in Conjecture~\ref{nw} would be best possible in the sense that for all
$n\ge 3$ and all $k<n/2$ there is a non-Hamiltonian strongly connected digraph~$G$ on~$n$ vertices
which satisfies the degree conditions except that $d^+_k,d^-_k\ge k+1$ are replaced by
$d^+_k,d^-_k\ge k$ in the $k$th pair of conditions.
To see this, take an independent set~$I$ of size $k<n/2$ 
and a complete digraph~$K$ of order~$n-k$. Pick a set~$X$ of~$k$ vertices of~$K$
and add all possible edges (in both directions) between~$I$ and~$X$. The digraph~$G$
thus obtained is strongly connected, not Hamiltonian and
$$\underbrace{k, \dots ,k}_{k \text{ times}}, \underbrace{n-1-k, \dots , n-1-k}_{n-2k
\text{ times}}, \underbrace{n-1, \dots , n-1}_{k \text{ times}}$$ is both the out- and
indegree sequence of~$G$. In contrast to the undirected case there exist examples with 
a similar degree sequence to the above but whose structure is quite different (see~\cite{KOTchvatal} and~\cite{CKKOsemi}).
This is one of the reasons which makes the directed case much harder than the undirected one.
In~\cite{CKKOsemi}, the following approximate version of Conjecture~\ref{nw} for large digraphs was proved. 
\begin{theorem}[Christofides, Keevash, K\"uhn and Osthus~\cite{CKKOsemi}] \label{CKKOchvatal}
For every $\beta > 0$ there exists an integer $n_0 = n_0(\beta)$
such that the following holds. Suppose that $G$ is a digraph on $n \geq
n_0$ vertices such that for all $i < n/2$
\begin{itemize}
\item[(i)] $d_i^+ \geq \min{ \{i + \beta n,n/2\} }$ or $d^-_{n-i - \beta n} \geq n-i$;
\item[(ii)] $d_i^- \geq \min{ \{i + \beta n, n/2 \} }$ or $d^+_{n-i - \beta n} \geq n-i$.
\end{itemize}
Then $G$ contains a Hamilton cycle.
\end{theorem}
This improved a recent result in~\cite{KOTchvatal}, where the degrees in the
first parts of these conditions were not `capped' at $n/2$.
The earlier result in~\cite{KOTchvatal} was derived from a result in~\cite{KKOexact} on the existence of a 
Hamilton cycle in an oriented graph satisfying a certain expansion property.
Capping the degrees at $n/2$ makes the proof far more difficult:
the conditions of Theorem~\ref{CKKOchvatal} only imply a  rather weak expansion property and 
there are many types of digraphs which almost satisfy the conditions but are not Hamiltonian.

The following weakening of Conjecture~\ref{nw} was posed earlier by Nash-Williams~\cite{ch2}.
It would yield a digraph analogue of P\'osa's theorem.
\begin{conj}[Nash-Williams~\cite{ch2}]\label{nw2}
Let $G$ be a digraph on $n \geq 3$ vertices such that $d^+ _i,d^-_i \geq i+1 $
for all $i <(n-1)/2$ and such that additionally
$d^+_{\lceil n/2\rceil},d^-_{\lceil n/2\rceil} \geq \lceil n/2\rceil$ when~$n$ is odd. 
Then~$G$ contains a Hamilton cycle.
\end{conj}
The previous example shows the degree condition would be best possible in the same sense as described there. 
The assumption of strong connectivity is not necessary in Conjecture~\ref{nw2},
as it follows from the degree conditions.
Theorem~\ref{CKKOchvatal} immediately implies a corresponding approximate version of Conjecture~\ref{nw2}.
In particular, for half of the vertex degrees (namely those whose value is $n/2$), 
the result matches the conjectured value.

\subsection{Chv\'atal-Erd\H{o}s type conditions}

Another sufficient condition for Hamiltonicity in undirected graphs which is just as fundamental
as those listed in the introduction is the Chv\'{a}tal-Erd\H{o}s theorem~\cite{chvatalerdos}:
suppose that $G$ is an undirected graph
with $n\geq 3$ vertices, for which the vertex-connectivity number $\kappa(G)$
and the independence  number $\alpha(G)$ satisfy
$\kappa(G)\geq \alpha(G)$, then $G$ has a Hamilton cycle.
Currently, there is no digraph analogue of this.
Given a digraph $G$, let $\alpha_0(G)$ denote 
the size of the largest set $S$ so that $S$ induces no edge
and let $\alpha_2(G)$ be
the size of the largest set $S$ so that $S$ induces no cycle of length 2.
So $\alpha_0(G) \le \alpha_2(G)$. 
$\alpha_0(G)$ is probably the more natural extension of the independence number to digraphs.
However, even the following basic question (already discussed e.g.~in~\cite{jacksonordaz}) is 
still open.
\begin{question}
Is there a function $f_0(k)$ so that every digraph with $\kappa(G) \ge f_0(k)$ and
$\alpha_0(G) \le k$ contains a Hamilton cycle?
\end{question}
Here the connectivity $\kappa(G)$ of a digraph is defined to be the size of the smallest
set of vertices $S$ so that $G-S$ is either not strongly connected or consists  of a single vertex. 
The following result shows that the analogous function for $\alpha_2(G)$ does exist.
\begin{theorem}[Jackson~\cite{jacksonchvatal}] \label{jacksonthm}
If $G$ is a digraph with
$$\kappa(G)\geq 2^{\alpha_{2}(G)}(\alpha_{2}(G)+2)!,$$
then $G$ has a Hamilton cycle.
\end{theorem}  
The proof involves a `reduction' of the problem to the undirected case.
As observed by Thomassen and Chakroun (see~\cite{jacksonordaz} again), 
there are non-Hamiltonian digraphs with $\kappa(G)= \alpha_2(G)=2$ and 
$\kappa(G)= \alpha_2(G)=3$. But it could well be that every digraph satisfying
$\kappa(G)\ge \alpha_2(G) \ge 4$ has a Hamilton cycle.
Even the following weaker conjecture is still wide open.
\begin{conj}[Jackson and Ordaz~\cite{jacksonordaz}] \label{hampath}
If $G$ is a digraph with $\kappa(G)\geq \alpha_{2}(G)+1$, then $G$ contains a Hamilton cycle.
\end{conj}
(In fact, they even conjectured that $G$ as above is pancyclic.)
Since the problem seems very difficult, even (say) a bound on $\kappa$ which is
polynomial in $\alpha_{2}$
in Theorem~\ref{jacksonthm}  would be interesting.


\section{Hamilton cycles in oriented graphs}

Recall that an \emph{oriented graph} is a directed graph with no $2$-cycles.
Results on oriented graphs seem even more difficult to obtain than results for the digraph case
(the Caccetta-H\"aggkvist conjecture on the girth of oriented graphs of large minimum outdegree
is a notorious example of this kind). 
In particular, most problems regarding Hamiltonicity of such graphs were open until recently
and many open questions still remain.

\subsection{Minimum degree conditions}

Thomassen~\cite{thomassen_79}
raised the natural question of determining the minimum
semidegree that forces a Hamilton cycle in
an oriented graph. Thomassen initially believed that the correct minimum semidegree bound should be $n/3$
(this bound is obtained by considering a `blow-up' of an oriented triangle).
However, H\"aggkvist~\cite{HaggkvistHamilton} later gave a construction which 
gives a lower bound of $\lceil (3n-4)/8 \rceil -1$: For $n$ of the form $n=4m+3$ where $m$ is odd,
we construct~$G$ on $n$ vertices as in Figure~3.
\begin{figure}\label{extremal2}
\centering\footnotesize
\includegraphics[scale=0.45]{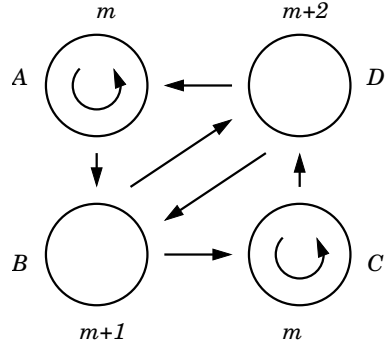}
\caption{An extremal example for Theorem~\ref{main}:
Partition the vertices into $4$ parts $A,B,C,D$, with $|A|=|C|=m$, $|B|=m+1$ and $|D|=m+2$.
Each of $A$ and $C$ spans a regular tournament, $B$ and $D$ are joined by a bipartite tournament
(i.e.~an orientation of the complete bipartite graph) which is as regular as possible.
We also add all edges from $A$ to $B$, from $B$ to $C$, from $C$ to $D$ and from $D$ to $A$.}
\end{figure}
Since every path which joins two vertices in~$D$ has to pass through~$B$, it follows that
every cycle contains at least as many vertices from~$B$ as it contains from~$D$.
As $|D|>|B|$ this means that one cannot cover all the vertices of~$G$ by disjoint cycles.
This construction can be extended to arbitrary~$n$
(see~\cite{KKOexact}). The following result exactly matches this bound and improves earlier ones of
several authors, e.g.~\cite{HaggkvistHamilton,HaggkvistThomasonHamilton,thomassenconj}. 
In particular, the proof builds on an approximate 
version which was proved in~\cite{kellyKO}.
\begin{theorem}[Keevash, K\"uhn and Osthus~\cite{KKOexact}] \label{main}
There exists an integer $n_0$ so that any oriented graph
$G$  on $n \ge n_0$ vertices with minimum
semidegree $\delta^0(G) \ge \frac{3n-4}{8}$ contains a Hamilton cycle.
\end{theorem}
Jackson conjectured that for regular oriented graphs one can significantly reduce the degree
condition.
\begin{conj}[Jackson~\cite{jacksonconj}] \label{jacksonconj}
For each $d>2$, every $d$-regular oriented graph $G$ on $n \le 4d+1$ vertices 
has a Hamilton cycle.
\end{conj}
The disjoint union of two regular tournaments on $n/2$ vertices shows that this 
would be best possible. Note that the degree condition is smaller than the one 
in Conjecture~\ref{diregconj}.
We believe that it may actually be possible to reduce the degree condition even further
if we assume that $G$ is strongly $2$-connected: is it true that for each $d>2$, every $d$-regular
strongly $2$-connected oriented graph $G$ on $n \le 6d$ vertices 
has a Hamilton cycle? A suitable orientation of the example in Figure~1 shows that this would 
be best possible.

\subsection{Ore-type conditions}
H\"aggkvist~\cite{HaggkvistHamilton} also made the following conjecture which is closely related
to Theorem~\ref{main}. Given an oriented graph~$G$, let~$\delta(G)$ denote the minimum degree of~$G$
(i.e.~the minimum number of edges incident to a vertex) and set
$\delta^*(G):=\delta(G)+\delta^+(G)+\delta^-(G)$.
\begin{conj}[H\"aggkvist~\cite{HaggkvistHamilton}] \label{haggconj}
Every oriented graph~$G$ on $n$ vertices with $\delta^*(G)>(3n-3)/2$ contains a Hamilton cycle.
\end{conj}
(Note that this conjecture does not quite imply Theorem~\ref{main} as it results in a marginally
greater minimum semidegree condition.)
In~\cite{kellyKO}, Conjecture~\ref{haggconj} was verified approximately, i.e.~if
$\delta^*(G) \ge (3/2+o(1))n$, then $G$ has a Hamilton cycle (note this implies an
approximate version of Theorem~\ref{main}).
The same methods also yield an approximate version of Ore's theorem for oriented graphs.
\begin{theorem}[Kelly, K\"uhn and Osthus~\cite{kellyKO}]\label{thm:Ore}
For every $\alpha>0$ there exists an integer $n_0=n_0(\alpha)$ such that every oriented
graph~$G$ of order $n\geq n_0$ with $d^+(x)+d^-(y)\ge (3/4+\alpha)n$ whenever $G$ does
not contain an edge from~$x$ to~$y$ contains a Hamilton cycle.
\end{theorem}
The construction in Figure~3 shows that the bound is best possible up to the term 
$\alpha n$. It would be interesting to obtain an exact version of this result.
\COMMENT{what extremal examples are there?}

Song~\cite{Songpancyclic} proved that every oriented graph on $n \ge 9$ vertices with $\delta(G) \ge n-2$
and $d^+(x)+d^-(y)\ge n -3$ whenever $G$ does
not contain an edge from~$x$ to~$y$ is pancyclic (i.e.~$G$ contains cycles of all possible lengths).
In~\cite{Songpancyclic} he also claims (without proof) that the condition is best possible for infinitely many $n$ as $G$ may fail to 
contain a Hamilton cycle otherwise. 
Note that Theorem~\ref{thm:Ore} implies that this claim is false.

\subsection{Degree sequence conditions and Chv\'atal-Erd\H{o}s type conditions}

In~\cite{KOTchvatal} a construction was described which showed that there is no satisfactory analogue
of P\'osa's theorem for oriented graphs: as soon as we allow a few vertices to have 
a degree somewhat below $3n/8$, then one cannot guarantee a Hamilton cycle.
The question of exactly determining all those degree sequences which guarantee a 
Hamilton cycle remains open though.

It is also not clear whether there may be a version of the Chv\'atal-Erd\H{o}s theorem 
for oriented graphs.


\section{Tournaments} \label{tournament}

A \emph{tournament} is an orientation of a complete graph. 
It has long been known that tournaments enjoy particularly strong Hamiltonicity properties: 
Camion~\cite{camion} showed that we only need to assume strong connectivity to ensure
that a tournament has a Hamilton cycle. 
Moon~\cite{moon} strengthened this by proving that every strongly connected tournament is even pancyclic.
It is easy to see that a minimum semidegree of $n/4$ forces a tournament on $n$ vertices to be strongly connected,
leading to a better degree condition for Hamiltonicity than that of $(3n-4)/8$ for the class of all oriented graphs.

\subsection{Edge-disjoint Hamilton cycles and decompositions}

A \emph{Hamilton decomposition} of a graph or digraph $G$ is a set of edge-disjoint Hamilton cycles which together
cover all the edges of~$G$. Not many examples of graphs with such decompositions are known.
One can construct a Hamilton decomposition of a complete graph if and only if its order is odd
(this was first observed by Walecki in the late 19th century).
Tillson~\cite{till} proved 
that a complete digraph $G$ on $n$ vertices has a Hamilton decomposition if and only if $n \not = 4,6$.  
The following conjecture of Kelly from 1968 (see Moon~\cite{moon}) would be a far-reaching
generalization of Walecki's result:
\begin{conj}[Kelly] \label{kelly}
Every regular tournament on $n$ vertices can be decomposed
into $(n-1)/2$ edge-disjoint Hamilton cycles.
\end{conj}
In~\cite{KOTkelly} we proved an approximate version of Kelly's conjecture. Moreover, the result 
holds even for oriented graphs $G$ which are not quite regular and whose `underlying' undirected
graph is not quite complete.
\begin{theorem}[K\"uhn, Osthus and Treglown~\cite{KOTkelly}] \label{kellythm}
For every $\eta_1 >0$ there exist $n_0= n_0 (\eta_1)$ and $\eta _2=\eta _2 (\eta_1) >0$ 
such that the following holds. Suppose that $G$ is an oriented graph on $n\geq n_0$ vertices such that
$\delta^0(G) \ge (1/2-\eta _2) n$.
Then $G$ contains at least $(1/2-\eta_1)n$ edge-disjoint Hamilton cycles. 
\end{theorem}
We also proved that the condition on the minimum semidegree can be relaxed to 
$\delta^0(G) \ge (3/8+\eta _2) n$.
This is asymptotically best possible since the construction described in Figure~3 is almost regular.

Some earlier support for Kelly's conjecture was provided by Thomassen~\cite{thom2}, who showed that the 
edges of every regular tournament can be covered by at most $12 n$ Hamilton cycles.
In this paper, we improve this to an asymptotically best possible result. 
We will give a proof (which relies on Theorem~\ref{kellythm}) in Section~\ref{coversec}.
\begin{theorem} \label{cover}
For every $\xi >0$ there exists an integer $n_0=n_0(\xi)$ 
so that every regular tournament $G$ on $n \geq n_0$ vertices contains a set of $(1/2+\xi)n$ Hamilton
cycles which together cover all the edges of $G$.
\end{theorem}
Kelly's conjecture has been generalized in several ways, 
e.g.~Bang-Jensen and Yeo~\cite{bangyeo} conjectured that every $k$-edge-connected tournament 
has a decomposition into $k$ spanning strong digraphs.
A bipartite version of Kelly's conjecture was also formulated by Jackson~\cite{jacksonconj}. 
Thomassen made the following conjecture which replaces the assumption of regularity by
high connectivity. 
\begin{conj}[Thomassen~\cite{thomassenconj}]
For every $k \ge 2$ there is an integer $f(k)$ so that every strongly $f(k)$-connected
tournament has $k$ edge-disjoint Hamilton cycles.
\end{conj}
A conjecture of Erd\H{o}s (see~\cite{thomassenconj}) which is also related to 
Kelly's conjecture states that almost all tournaments $G$
have at least $\delta^0(G)$ edge-disjoint Hamilton cycles.

Similar techniques as in the proof of the approximate version of Kelly's conjecture
were used at the same time in~\cite{cko} to prove approximate versions of two long-standing
conjectures of Nash-Williams on edge-disjoint Hamilton cycles in (undirected) graphs.
One of these results states that one can almost decompose any dense regular graph into Hamilton cycles.
\begin{theorem}[Christofides, K\"uhn and Osthus~\cite{cko}]\label{Regular}
For every $\eta > 0$ there is an integer $n_0=n_0(\eta)$ so that every
$d$-regular graph on $n \geq n_0$ vertices with $d \geq (1/2 +
\eta)n$ contains at least $(d - \eta n)/2$ edge-disjoint
Hamilton cycles.
\end{theorem}
In Section~\ref{coversec} we deduce the following analogue of Theorem~\ref{cover}:
\begin{theorem} \label{RegularII}
For every $\xi > 0$ there is an integer $n_0=n_0(\xi)$ so that every
$d$-regular graph $G$ on $n \geq n_0$ vertices with $d \geq (1/2 +
\xi)n$ contains at most  $(d + \xi n)/2$ 
Hamilton cycles which together cover all the edges of $G$.
\end{theorem}


\subsection{Counting Hamilton cycles in tournaments} 

One of the earliest results on tournaments (and the probablistic method), was obtained by
Szele~\cite{szele}, who showed that the maximum number $P(n)$ of Hamilton paths in a tournament on $n$ vertices 
satisfies $P(n)=O(n!/2^{3n/4})$ and $P(n) \ge n!/2^{n-1}=:f(n)$.
The lower bound is obtained by considering a random tournament.
The best upper bound is due to Friedgut and Kahn~\cite{friedgutkahn}
who showed that $P(n) = O(n^{c}f(n))$, where $c$ is slightly less than $5/4$.
The best current lower bound is due to Wormald~\cite{wormald}, who showed that
$P(n) \ge (2.855+o(1))f(n)$. So in particular, $P(n)$ is not attained for random tournaments.
Also, he conjectured that this bound is very close to the correct value.


Similarly, one can define the maximum number $C(n)$ of Hamilton cycles in a tournament on
$n$ vertices. 
Note that by considering a random tournament again, we obtain $C(n) \ge (n-1)!/2^n=:g(n)$.
Unsurprisingly, $C(n)$ and $P(n)$ are very closely related,
e.g.~we have $P(n) \ge n C(n)$.
In particular, the main result in~\cite{friedgutkahn} states that $C(n) = O(n^{c}g(n))$, where
$c$ is the same as above. This implies the above bound on $P(n)$, since 
Alon~\cite{alontournament} observed that $P(n) \le 4C(n+1)$.
Also, Wormald~\cite{wormald} showed that $C(n) \ge (2.855+o(1))g(n)$.
(Note this also follows by combining Alon's observation with the lower bound on 
$P(n)$ in~\cite{wormald}.)

Of course, in general it does not make sense to ask for the minimum number of 
Hamilton paths or cycles in a tournament.
However, the question does  make sense for regular tournaments.
Friedgut and Kahn~\cite{friedgutkahn} asked whether 
the number of Hamilton cycles in a regular tournament is always at least
$\Omega(g(n))$.
The best result towards this was recently obtained by Cuckler~\cite{cuckler}, who
showed that every regular tournament on $n$ vertices contains at least
$n!/(2+o(1))^n$ Hamilton cycles. This also answers an earlier question of Thomassen.
Asking for the minimum number of Hamilton paths in a tournament $T$ also
makes sense if we assume that $T$ is strongly connected.
Busch~\cite{busch} determined this number exactly by showing that an earlier construction of
Moon is best possible.
The related question on the minimum number of Hamilton cycles in a strongly $2$-connected
tournament is still open (see~\cite{busch}).

\subsection{Sumner's universal tournament conjecture}

Sumner's universal tournament conjecture states that every tournament on $2n-2$ vertices
contains every tree on $n$ vertices.
In~\cite{KMOsumner} an approximate version of this conjecture was proved
and subsequently in~\cite{exactsumner}, the conjecture was proved for all large trees
(see e.g.~\cite{KMOsumner} for a discussion of numerous previous results).
The proof in~\cite{exactsumner} builds on several structural results proved in~\cite{KMOsumner}.
\begin{theorem}[K\"uhn, Mycroft and Osthus~\cite{KMOsumner, exactsumner}] \label{sumnerapprox}
There is an integer $n_0$ such that for all $n \geq n_0$ every
tournament $G$ on $2n-2$ vertices contains any directed tree $T$ on $n$ vertices.
\end{theorem}
While this result is not directly related to the main topic of the survey (i.e.~Hamilton cycles), there are several connections.
Firstly, just as with many of the new results in the other sections, the concept 
of a robust expander is crucial in the proof of Theorem~\ref{sumnerapprox}.
Secondly, the proof of Theorem~\ref{sumnerapprox} also makes direct use of the fact that a robust expander contains a 
Hamilton cycle (Theorem~\ref{expanderthm}). Suitable parts of the tree $T$ are
embedded by considering a random walk on (the blow-up of) such a Hamilton cycle.

In \cite{KMOsumner}, we also proved that if $T$ has bounded maximum degree, then it suffices if the tournament $G$ has $(1+\alpha)n$
vertices. This is best possible in the sense that the `error term' $\alpha n$ cannot be completely omitted
in general.
But it seems possible that it can be reduced to a constant which depends only on the maximum degree of $T$.
If $T$ is an orientation of a path, then the error term can be omitted completely:
Havet and Thomass\'e~\cite{HTrosenfeld} proved that every tournament
on at least 8 vertices contains every possible orientation of a Hamilton path
(for arbitrary orientations of Hamilton cycles, see Section~\ref{arbitraryorient}).

\section{Generalizations} \label{general}

In this section, we discuss several natural ways of strengthening the notion of a Hamilton cycle.

\subsection{Pancyclicity} 

Recall that a graph (or digraph) is \emph{pancyclic} if it contains a cycle of every possible length.
Dirac's theorem implies that a graph on $n \ge 3$ vertices is pancyclic if it has 
minimum degree greater than $n/2$. (To see this, remove a vertex $x$ and apply Dirac's theorem to the remaining
subgraph to obtain a cycle of length $n-1$. Then consider the neighbourhood of $x$ on this cycle to obtain cycles
of all possible lengths through~$x$.)
Similarly, one can use Ghouila-Houri's theorem to deduce that every digraph on $n$ vertices with  minimum
semidegree greater than $n/2$ is pancyclic.
In both cases, the complete bipartite (di-)graph whose vertex class sizes
are as equal as possible shows that the bound is best possible.
More generally, the same trick also works for Meyniel's theorem:
let $G$ be a strongly connected digraph on $n \ge 2$ vertices.
If $d(x)+d(y) \ge 2n+1$ for all pairs of non-adjacent vertices $x \neq y$
in $G$, then $G$ is pancyclic. (Indeed, the conditions imply that either $G$
contains a strongly connected tournament or contains a vertex $x$ with $d(x)>n$, 
in which case we can proceed as above.)%
\COMMENT{$d_{G'}(x) \ge n+ 1/2$ and $d(x)+d(y) \ge 2n+1-4=2(n-1)-1$ in $G'$ and $x$ has min in- and outdegree $\ge 1$}
If $n$ is even, the bound $2n+1$ is best possible. For $n$ is odd, it follows from a result of 
Thomassen~\cite{tom} that one can improve it to $2n$.


For oriented graphs the minimum semidegree threshold which guarantees pancyclicity
turns out to be $(3n-4)/8$, i.e.~the same threshold as for Hamiltonicity (see~\cite{kellyKOpan}).
The above trick of removing a vertex does not work here. 
Instead, to obtain `long' cycles one can modify the proof of Theorem~\ref{main}.
A triangle is guaranteed by results on the Caccetta-H\"aggkvist conjecture --
e.g.~a very recent result of Hladk\'y, Kr\'al and Norine~\cite{HKN} states that every oriented graph on $n$ vertices with 
minimum semidegree at least $0.347n$ contains a 3-cycle.
Short cycles of length $\ell\ge 4$ can be guaranteed by a result in~\cite{kellyKOpan} which states that
for all $n\ge 10^{10}\ell$ every oriented graph~$G$ on $n$ vertices with
$\delta^0(G)\ge \lfloor n/3\rfloor+1$
contains an $\ell$-cycle. This is best possible for all those $\ell\geq 4$ which are not
divisible by~3.
Surprisingly, for some other values of~$\ell$, an $\ell$-cycle is forced by a much
weaker minimum degree condition. In particular, the following conjecture was made in~\cite{kellyKOpan}.
\begin{conj}[Kelly, K\"uhn and Osthus~\cite{kellyKOpan}] \label{con:short}
Let $\ell\geq 4$ be a positive integer and let~$k$ be the
smallest integer that is greater than~$2$ and does not
divide~$\ell$.
Then there exists an integer $n_0=n_0(\ell)$ such that every oriented
graph~$G$ on $n\geq n_0$ vertices with minimum semidegree $\delta^0(G)\geq \lfloor n/k\rfloor+1$
contains an $\ell$-cycle.
\end{conj}
The extremal examples for this conjecture are always `blow-ups' of cycles of length $k$.
Possibly one can even weaken the condition by requiring only the outdegree of $G$ to be large.
It is easy to see that
the only values of $k$ that can appear in Conjecture~\ref{con:short} are of the form
$k=p^s$ with $k \ge 3$, where $p$ is a prime and $s$ a positive integer.%

\subsection{Arbitrary orientations} \label{arbitraryorient}

As mentioned earlier, the most natural notion of a cycle in a digraph is to have all edges 
directed consistently. But it also makes sense to ask for Hamilton cycles where the edges
are oriented in some prescribed way, e.g.~to ask for an `antidirected' Hamilton cycle
where consecutive edges have opposite directions.
Surprisingly, it turns out that both for digraphs and oriented graphs the minimum 
degree threshold which guarantees a `consistent' Hamilton cycle is approximately the 
same as that which guarantees an arbitrary orientation of a Hamilton cycle.
\begin{theorem}[H\"aggkvist and Thomason~\cite{HTdigraph}]
There exists an $n_0$ so that every digraph $G$ on $n\ge n_0$ vertices with
minimum semidegree $\delta^0(G) \ge n/2+n^{5/6}$ contains every orientation of a Hamilton cycle.
\end{theorem}
In~\cite{HaggkvistThomasonHamilton},
they conjectured an analogue of this for oriented graphs, which was recently proved by Kelly.
\begin{theorem}[Kelly~\cite{Kelly}]
For every $\alpha>0$ there exists an integer $n_0=n_0(\alpha)$ such that every oriented
graph~$G$ on $n\geq n_0$ vertices with minimum semidegree $\delta^0(G)\ge (3/8+\alpha)n$
contains every orientation of a Hamilton cycle.
\end{theorem}
The proof of this result uses Theorem~\ref{main} as the well as the notion of expanding digraphs.
Interestingly, Kelly observed that the thresholds for various orientations do not coincide exactly:
for instance, if we modify the example in Figure~3 so that all classes have the same odd size, then 
the resulting oriented graph has minimum semidegree $(3n-4)/8$ but no antidirected Hamilton cycle.

Thomason~\cite{thomasontournament} showed that for large tournaments strong
connectivity ensures every possible orientation
of a Hamilton cycle. More precisely, he showed that for $n \ge 2^{128}$, every tournament on $n$
vertices contains all possible orientations of a Hamilton cycle, except possibly the `consistently oriented' one.
(Note that this also implies that every large tournament contains every orientation of a Hamilton path,
i.e.~a weaker version of the result in~\cite{HTrosenfeld} mentioned earlier.)
The bound on $n$ was later reduced to 68 by Havet~\cite{havettournament}.
Thomason conjectured that the correct bound is $n\ge 9$.

\subsection{$k$-ordered Hamilton cycles}


Suppose that we require our (Hamilton) cycle to visit several vertices in a specific order. 
More formally, we say that a graph $G$ is \emph{$k$-ordered} if 
for every sequence $s_1,\dots,s_k$ of distinct vertices of~$G$ there is a cycle
which encounters $s_1,\dots,s_k$ in this order. $G$~is \emph{$k$-ordered Hamiltonian}
if it contains a Hamilton cycle with this property. 
Kierstead, S\'ark\"ozy and Selkow~\cite{KSSordered} determined the minimum 
degree which forces an (undirected) graph to be $k$-ordered Hamiltonian.
\begin{theorem}[Kierstead, S\'ark\"ozy and Selkow~\cite{KSSordered}] \label{kordered}
For all $k \ge 2$, 
every graph on $n \ge 11k-3$ vertices of minimum degree at least $\lceil n/2 \rceil + \lfloor k/2 \rfloor-1$
is $k$-ordered Hamiltonian. 
\end{theorem}
The extremal example consists of two cliques intersecting in $k-1$ vertices if $k$ is even 
and two cliques intersecting in $k-2$ vertices if $k$ is odd.
The case when $n$ is not too large compared to $k$ is still open. The corresponding 
Ore-type problem was solved in~\cite{orekordered}.
Here the Ore-type result does not imply the Dirac-type result above.
Many variations and stronger notions have been investigated (see~e.g.~\cite{gould} again).

Directed graphs form a particularly natural setting for this kind of question.
The following result gives a directed analogue of Theorem~\ref{kordered}.
\begin{theorem}[K\"uhn, Osthus and Young~\cite{KOY}]\label{thm:ordHam}
For every $k \ge 3$ there is an integer $n_0=n_0(k)$ such that every 
digraph~$G$ on~$n\ge n_0$ vertices with minimum semidegree $\delta^0(G)\ge \lceil(n+k)/2\rceil -1$
is $k$-ordered Hamiltonian.
\end{theorem}
Note that if~$n$ is even and~$k$ is odd the bound on the minimum semidegree is
slightly larger than in the undirected case.
However, it is best possible in all cases. In fact, if the minimum semidegree is smaller, 
it turns out that $G$ need not even be $k$-ordered. 
Again, the family of extremal examples turns out to be much richer than in the undirected case.
Note that every Hamiltonian digraph is $2$-ordered Hamiltonian, so the case when $k \le 2$ in Theorem~\ref{thm:ordHam}
is covered by Ghouila-Houri's theorem.
It would be interesting to obtain an Ore-type or an oriented version of Theorem~\ref{thm:ordHam}.%

\subsection{Factors with prescribed cycle lengths}

Another natural way of generalizing Dirac's theorem is to ask for a certain set of vertex-disjoint cycles in $G$
which together cover all the vertices of $G$ (note this also
generalizes the notion of pancyclicity). For large undirected graphs, Abassi~\cite{Abbasi} determined the 
minimum degree which guarantees $k$ vertex-disjoint cycles in a graph $G$ whose (given) lengths are $n_1,\dots,n_k$,
where the $n_i$ sum up to $n$ and where the order $n$ of $G$ is sufficiently large.
As in the case of Hamilton cycles, the corresponding questions for directed and
oriented graphs appear much harder than in the undirected case 
and again much less is known. 
Keevash and Sudakov~\cite{keevashsudakov} recently obtained the following result.
\begin{theorem}[Keevash and Sudakov~\cite{keevashsudakov}]
There exist positive constants $c,C$ and an integer~$n_0$ so that whenever $G$ is an
oriented graph on $n\ge n_0$ vertices
with minimum semidegree at least $(1/2-c)n$ and whenever $n_1,\dots,n_t$ are so that 
$\sum_{i=1}^t n_i \le n-C$, then $G$ contains vertex-disjoint cycles of length $n_1,\dots,n_t$.
\end{theorem}
In general, one cannot take $C=0$.
In the case of triangles (i.e.~when all the $n_i=3$), they show that one can choose $C=3$.
This comes very close to proving a recent conjecture formulated independently by Cuckler and  
Yuster~\cite{yuster}, which states that 
every regular tournament on $n=6k+3$ vertices contains vertex-disjoint triangles covering all the vertices
of the tournament.
Similar questions were also raised earlier by Song~\cite{song}. For instance,
given $t$, he asked for the smallest integer $f(t)$ so that all but a finite number of strongly $f(t)$-connected
tournaments $T$ satisfy the following: Let $n$ be the number of vertices of $T$ and let 
$\sum_{i=1}^t n_i = n$. Then $T$ contains vertex-disjoint cycles of length $n_1,\dots,n_t$.
Chen, Gould and Li~\cite{cgl} proved the weaker result that every sufficiently large
$t$-connected tournament $G$ contains $t$ vertex-disjoint cycles which together cover all the vertices of $G$.
This proved a conjecture of Bollob\'as.

\subsection{Powers of Hamilton cycles}

Sark\"ozy, Koml\'os and Szemer\'edi~\cite{KSSz98}
showed that every sufficiently large graph $G$ on $n$ vertices with minimum degree at least
$kn/(k+1)$ contains the $k$th power of a Hamilton cycle.
Extremal examples are complete $(k+1)$-partite graphs with classes of almost equal size.
It appears likely that the situation for digraphs is similar.
However, just as for ordinary Hamilton cycles, it seems that for oriented graphs
the picture is rather different.
(Both for digraphs and oriented graphs, the most natural definition of the $k$th power of 
a cycle is a cyclically ordered set of vertices so that every vertex sends an edge to the next
$k$ vertices in the ordering.)
\begin{conj}[Treglown~\cite{treglownthesis}]\label{squareconj}
For every $\eps>0$ there is an integer $n_0=n_0(\eps)$ so that every oriented graph $G$
on $n \ge n_0$ vertices with minimum semidegree at least $(5/12 + \eps)n$ contains the square of a Hamilton cycle.
\end{conj}
A construction which shows that the constant $5/12$ cannot be improved is given in Figure~4.
\begin{figure}\label{square}
\centering\footnotesize
\includegraphics[scale=0.42]{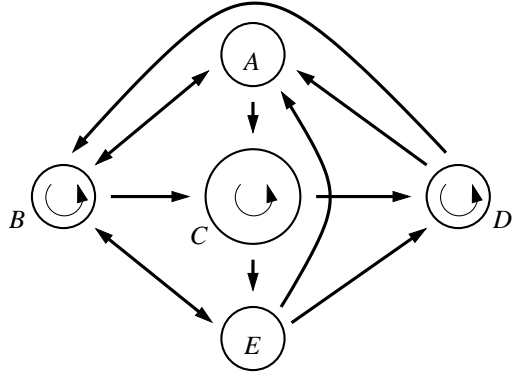}
\caption{An extremal example for Conjecture~\ref{squareconj}:
The set sizes are $|A|=m$, $|B|=m-1$, $|C|=2m+1$, $|D|=m-1$ and $|E|=m+1$, where $m$ is even.
$B$, $C$ and $D$ induce regular tournaments, while $A$ and $E$ induce independent sets. A single arrow 
(e.g.~from $B$ to $C$) indicates an orientation of the complete bipartite graph from 
the initial set towards the final set.
A double edge (e.g.~between $B$ and $E$) indicates an orientation of the complete bipartite graph
so that within each set, the in- and outdegrees of the vertices differ by at most one.}
\end{figure}
We claim that the square of any Hamilton cycle would have to visit a vertex of $B$ in between two visits of $E$.
Since $|B|<|E|$, this shows that the graph does not contain the square of a Hamilton cycle.
To prove the claim, suppose that $F$ is a squared Hamilton cycle and consider a vertex $e$ of $F$ which lies in $E$.
Then the predecessor of $e$ lies in $C$ or $B$, so without loss of generality we may assume 
that it is a vertex $c_1 \in C$. Again, the predecessor of $c_1$ lies in $C$ or $B$
(since it must lie in the common inneighbourhood of $C$ and $E$), 
so without loss of generality we may assume that it is a vertex $c_2 \in C$.
The predecessor of $c_2$ can now lie in $A$, $B$ or $C$. If it lies in $B$ we are done again,
if it is a vertex $c_3 \in C$, we consider its predecessor, which can again only lie in $A$, $B$ or $C$.
Since $F$ must visit all vertices, it follows that we eventually arrive at a predecessor $a \in A$
whose successor on $F$ is some vertex $c \in C$. The predecessor of $a$ on $F$ must lie in the common
inneighbourhood of $a$ and $c$, so it must lie in $B$, as required.

For the case of tournaments, the problem was solved asymptotically by Bollob\'as and H\"aggkvist~\cite{BHpower}.
Given a tournament $T$ of large order~$n$ with minimum semidegree at least
$n/4+\eps n$, they proved that (for fixed $k$) $T$  contains the 
$k$th power of a Hamilton cycle. So asymptotically, the semidegree threshold for an 
ordinary Hamilton cycle in a tournament is the same as that for the $k$th power of a Hamilton cycle.



\section{Robustly expanding digraphs} \label{expansion}

Roughly speaking, a graph is an expander if for every set $S$ of vertices the 
neighbourhood $N(S)$ of $S$ is significantly larger than $S$ itself. 
A number of papers have recently demonstrated that there is a remarkably close 
connection between Hamiltonicity and expansion (see e.g.~\cite{HKS}).
The following notion of robustly expanding (dense) digraphs was introduced in~\cite{KOTchvatal}.

Let $0<\nu\le  \tau<1$. Given any digraph~$G$ on $n$ vertices and
$S\subseteq V(G)$, the \emph{$\nu$-robust outneighbourhood~$RN^+_{\nu,G}(S)$ of~$S$}
is the set of all those vertices~$x$ of~$G$ which have at least $\nu n$ inneighbours
in~$S$. $G$ is called a \emph{robust $(\nu,\tau)$-outexpander}
if $|RN^+_{\nu,G}(S)|\ge |S|+\nu n$ for all
$S\subseteq V(G)$ with $\tau n< |S|< (1-\tau)n$.
As the name suggests, this notion has the advantage that it is preserved even if we
delete some vertices and edges from $G$.
We will also use the more traditional (and weaker) notion of a \emph{$(\nu,\tau)$-outexpander}, which means
$|N^+(S)|\ge |S|+\nu n$ for all $S\subseteq V(G)$ with $\tau n<|S|< (1-\tau)n$.
\begin{theorem}[K\"uhn, Osthus and Treglown~\cite{KOTchvatal}]\label{expanderthm}
Let $n_0$ be a positive integer and $\nu,\tau,\eta$ be
positive constants such that $1/n_0\ll\nu\le \tau\ll\eta<1$. Let~$G$ be a digraph on~$n\ge n_0$ vertices with
$\delta^0(G)\ge \eta n$ which is a robust $(\nu,\tau)$-outexpander. Then~$G$ contains a Hamilton cycle.
\end{theorem}
Theorem~\ref{expanderthm} is used in~\cite{KOTchvatal} to give a weaker version of Theorem~\ref{CKKOchvatal}
(i.e.~without the degrees capped at $n/2$). In the same paper it is also applied to prove a conjecture
of Thomassen regarding a weak version of Conjecture~\ref{kelly} (Kelly's conjecture).
One can also  use it to prove e.g.~Theorem~\ref{thm:Ore} and thus an approximate version of Theorem~\ref{main}. 
(Indeed, as proved in~\cite{kellyKO}, the degree conditions of Theorem~\ref{thm:Ore} imply expansion, the proof for 
robust expansion is similar.)
As mentioned earlier, it is also used as a tool in the proof of 
Theorem~\ref{sumnerapprox}. 
Finally, we will also use it in the next subsection to prove Theorem~\ref{cover}.

In~\cite{KOTchvatal}, Theorem~\ref{expanderthm} was deduced from a result in~\cite{KKOexact}.
The proof of the result in~\cite{KKOexact} (and a similar approach in~\cite{kellyKO}) 
in turn relied on Szemer\'edi's regularity lemma and a (rather technical)
version of the Blow-up lemma due to Csaba~\cite{Csababipartite}.
A (parallel) algorithmic version of Theorem~\ref{expanderthm} was also proved in~\cite{CKKOalgo}.
Below, we give a brief sketch of a proof of Theorem~\ref{expanderthm} which avoids any use of the Blow-up lemma
and is based on an approach in~\cite{CKKOsemi}.

The {\em density} of a bipartite graph $G$ with vertex
classes $A$ and $B$ is defined to be $d(A,B) =
\frac{e(A,B)}{|A||B|}$, where $e(A,B)$ denotes the number of edges between $A$ and $B$.
Given $\eps > 0$, we say that $G$ is $\eps$-{\em
regular} if for all subsets $X \subseteq  A$ and $Y \subseteq B$ with
$|X| \geq \eps|A|$ and $|Y| \geq \eps |B|$ we have that $|d(X,Y) -
d(A,B)| < \eps$. We also say that $G$ is $(\eps,d)$-{\em super-regular} if it is
$\eps$-regular and furthermore every vertex $a \in A$ has degree at least $d|B|$
and similarly for every $b \in B$. These definitions generalize naturally to non-bipartite (di-)graphs.

We also need the result that every super-regular digraph contains a Hamilton cycle.
%
\begin{lemma} \label{fk} 
Suppose that $1/n_0\ll\eps \ll d \ll 1$ and $G$ is
an $(\eps,d)$-super-regular digraph on $n\ge n_0$ vertices.
Then $G$ contains a Hamilton cycle.
\end{lemma}
Lemma~\ref{fk} is a special case e.g.~of a result of 
Frieze and Krivelevich~\cite{FK}, who proved that an $(\eps,d)$-super-regular
digraph on $n$ vertices has almost $dn$ edge-disjoint Hamilton cycles if $n$ is large. Here we also give a sketch of
a direct proof of Lemma~\ref{fk}.

We first prove that $G$ contains a 1-factor. Consider the auxiliary bipartite graph whose vertex
classes $A$ and $B$ are copies of $V(G)$
with an edge between $a\in A$ and $b \in B$ if there is an edge from $a$ to $b$ in $G$.
One can show that this bipartite graph has a perfect matching (by Hall's marriage theorem),
which in turn corresponds to a $1$-factor in $G$. 

It is now not hard to prove the lemma using the `rotation-extension' technique:
Choose a 1-factor of $G$.
Now remove an edge of a cycle in this 1-factor and let $P$ be the resulting path.
If the final vertex of $P$ has any outneighbours on another cycle $C$ of the 1-factor, we can extend $P$
into a longer path which includes the vertices of $C$
(and similarly for the initial vertex of $P$). 
We repeat this as long as possible (and one can always ensure that the extension step
can be carried out at least once). So we may assume that all outneighbours of the final vertex of $P$
lie on $P$ and similarly for the initial vertex of $P$. 
Together with the $\eps$-regularity this can be used to find a cycle with the same vertex set as $P$.
Eventually, we arrive at a Hamilton cycle.

\medskip

\noindent
{\bf Sketch proof of Theorem~\ref{expanderthm}}.
Choose $\eps,d$ to satisfy $1/n_0 \ll \eps \ll d \ll \nu$.
The first step is to apply a directed version of Szemer\'edi's regularity lemma to $G$.
This gives us a partition of the vertices of $G$ into clusters $V_1,\dots,V_k$ and
an exceptional set $V_0$ so that $|V_0| \le \eps n$ and all the clusters have size $m$.
Now define a `reduced' digraph $R$ whose vertices are the clusters $V_1,\dots,V_k$ and with 
an edge from $V_i$ to $V_j$ if the bipartite graph spanned by the edges from $V_i$ to $V_j$
is $\eps$-regular and has density at least $d$.
Then one can show (see Lemma~14 in~\cite{KOTchvatal}) that $R$ is still a $(\nu/2,2\tau)$-outexpander
(this is the point where we need the robustness of the expansion in $G$) with minimum semidegree at least
$\eta k/2$. This in turn can be used to show that $R$ has a $1$-factor $F$ (using the same
auxiliary bipartite graph as in the proof of Lemma~\ref{fk}).
By removing a small number of vertices from the clusters, we can also assume that the 
bipartite subgraphs spanned by successive clusters on each cycle of $F$ are super-regular, 
i.e.~have high minimum degree. For simplicity, assume that the cluster size is still $m$.

Moreover, 
since $G$ is an expander, we can find a short path in $G$ between clusters of different cycles 
of $F$ and also between any pair of exceptional vertices. 
However, we need  to choose such paths without affecting any of the useful structures
that we have found so far. For this, we will consider paths which `wind around' cycles in $F$
before moving to another cycle. More precisely,
a {\em shifted walk} from a cluster $A$ to a cluster $B$
is a walk $W(A,B)$ of the form
\[ W(A,B) = X_1 C_1 X^-_1 X_2 C_2 X^-_2 \dots X_t C_t X^-_t X_{t+1},\]
where $X_1=A$, $X_{t+1} = B$, $C_i$ is the  cycle of $F$ containing $X_i$,
and for each $1 \leq i \leq t$, $X^-_i$ is the predecessor of $X_i$ on $C_i$
and the edge $X^-_i X_{i+1}$ belongs to $R$.
We say that $W$ as above \emph{traverses $t$ cycles} (even if some $C_i$ appears several times in $W$).
We also say that the clusters $X_2,\dots,X_{t+1}$ are the \emph{entry clusters}
(as this is where $W$ `enters' a cycle $C_i$)
and the clusters $X_1^-,\dots,X_{t}^-$ are the \emph{exit clusters} of $W$.
Note that 
\begin{itemize}
\item[(i)] for any cycle of $F$, its clusters are visited the same number of times by $W(A,B)- B$.
\end{itemize}
Using the expansion of $R$, it is not hard to see that
\begin{itemize}
\item[(ii)]for any clusters $A$ and $B$ there is a 
shifted walk from $A$ to $B$ which does not traverse too many cycles.
\end{itemize}
Indeed, the expansion property implies that the number of clusters one can reach by 
traversing $t$ cycles is at least $t \nu k/2$ as long as this is significantly less than the 
total number $k$ of clusters.

Now we will `join up' the exceptional vertices using shifted walks. For this, write $V_0=\{ a_1,\dots, a_{\ell} \}$. 
For each exceptional vertex $a_i$ choose a cluster $T_i$ so that $a_i$ has many outneighbours in $T_i$. 
Similarly choose a cluster $U_i$ so that $a_i$ has many inneighbours in $U_i$ 
and so that
\begin{itemize}
\item[(iii)] no cluster appears too often as a $T_i$ or a $U_i$.
\end{itemize}
Given a cluster $X$, let $X^-$ be the predecessor of $X$ on the cycle of $F$ which contains $X$
and let $X^+$ be its successor.
Form a `walk' $W$ on $V_0 \cup V(R)$ which starts at $a_1$, then moves to $T_1$, then follows a shifted walk
from $T_1$ to $U^+_2$, then it winds around the entire cycle of $F$ containing $U^+_2$ until it reaches $U_2$.
Then $W$ moves to $a_2$, then to $a_3$ using a shifted walk as above until it has visited all the exceptional 
vertices. Proceeding similarly, we can ensure that $W$ has the following properties:
\begin{itemize}
\item[(a)] $W$ is a closed walk which visits all of $V_0$ and all of $V(R)$.
\item[(b)] For any cycle of $F$, its clusters are visited the same number of times by~$W$.
\item[(c)] Every cluster appears at most $m/10$ times as an entry or exit cluster.
\end{itemize} 
(b) follows from (i) and (c) follows from (ii) and~(iii).
The next step towards a Hamilton cycle would be to find a cycle $C$ in $G$ which corresponds to $W$
(i.e.~each occurrence of a cluster in $W$ is replaced by a distinct vertex of $G$ lying in this cluster).
Unfortunately, the fact that $V_0$ may be much larger than the cluster size $m$ implies that 
there may be clusters which are visited more than $m$ times by $W$, which makes it impossible to find 
such a $C$. So we will apply a `short-cutting' technique to~$W$ which avoids `winding around' the cycles of $F$
too often.

For this, we now  fix edges in $G$ corresponding to all those edges of $W$ that do not lie within a cycle of $F$.
These edges of $W$ are precisely the edges in $W$ at the exceptional vertices as well as
all the edges of the form $AB$ where $A$ is used as an exit cluster by $W$ and $B$ is used as an
entrance cluster by $W$.
For each edge $a_iT_i$ at an exceptional vertex we choose an edge $a_ix$, where $x$ is an
outneighbour of $a_i$ in $T_i$. We similarly choose an edge $ya_i$ from $U_i$ to~$a_i$ for each $U_ia_i$.
We do this in such a way that all these edges are disjoint outside~$V_0$. 
For each occurrence of $AB$ in $W$, where $A$ is used as an exit cluster by $W$ and $B$ is used as an
entrance cluster, we choose an edge $ab$ from  $A$ to $B$ in $G$ so that all these
edges are disjoint from each other and from the edges chosen for the exceptional vertices (we use~(c) here). 

Given a cluster $A$,  let $A_{entry}$ be the set of all those vertices in $A$ which
are the final vertex of an edge of $G$ fixed so far
and let $A_{exit}$ be the set of all those vertices in $A$ which are the initial vertex of an edge
of $G$ fixed so far. So $A_{entry} \cap A_{exit} = \emptyset$.
Let $G_A$ be the bipartite graph whose vertex classes are $A \setminus A_{exit}$ and
$A^+ \setminus A^+_{entry}$ and whose edges are all the edges from $A \setminus A_{exit}$
to $A^+ \setminus A^+_{entry}$ in $G$. Since $W$ consists of shifted walks, it is easy to see that the vertex
classes of $G_A$ have equal size. Moreover, it is possible to carry out the previous steps in such a way that
$G_A$ is super-regular (here we use~(c) again). This in turn means that $G_A$ has a perfect matching $M_A$.
These perfect matchings (for all clusters $A$) together with all the edges of $G$ fixed so far form a
$1$-factor ${\mathcal C}$ of $G$.
It remains to transform ${\mathcal C}$ into a Hamilton cycle.

We claim that for any cluster $A$, we can find a perfect matching $M'_A$ in $G_A$
so that if we replace $M_A$ in ${\mathcal C}$ with $M'_A$, then all
vertices of $G_A$ will lie on a common cycle in the new
$1$-factor~${\mathcal C}$.

To prove this claim we proceed as follows. For every $a \in A^+
\setminus  A^+_{entry}$, we move along the cycle $C_a$ of ${\mathcal C}$ containing $a$
(starting at $a$) and let $f(a)$ be the first vertex on $C_a$ in
$A  \setminus A_{exit}$. Define an auxiliary digraph $J$ on $A^+\setminus A^+_{entry}$ such that 
$N^+_J(a):=N^+_{G_A}(f(a))$.
So $J$ is obtained by identifying each pair
$(a,f(a))$ into one vertex with an edge from $(a,f(a))$ to
$(b,f(b))$ if $G_A$ has an edge from $f(a)$ to $b$. Since $G_A$ is
super-regular, it follows that $J$ is also super-regular.  By Lemma~\ref{fk}, $J$ has a Hamilton cycle,
which clearly corresponds to a perfect matching $M'_A$ in~$G_A$ with
the desired property.

We now repeatedly apply the above claim to every cluster. Since $A_{entry} \cap A_{exit} = \emptyset$
for each cluster $A$, this ensures that all vertices which lie in clusters on the same cycle of $F$ will lie
on the same cycle of the new $1$-factor ${\mathcal C}$. Since by~(a) $W$ visits all clusters, this in turn implies
that all the non-exceptional vertices will lie in the same cycle of  ${\mathcal C}$.
Since the exceptional vertices form an
independent set in~$\mathcal{C}$, it follows that $\mathcal{C}$ is actually a 
Hamilton cycle.
\endproof

\subsection{Covering regular graphs and tournaments with Hamilton cycles} \label{coversec}

Here we give proofs of Theorems~\ref{cover} and~\ref{RegularII}.
The proof of Theorem~\ref{cover} uses Theorems~\ref{kellythm} and~\ref{expanderthm}.

\medskip

\noindent{\bf Proof of Theorem~\ref{cover}.} Choose new constants $\eta_1,\nu,\tau$ such that
$1/n_0\ll \eta_1\ll \nu\le \tau\ll \xi$. Consider any regular tournament $G$ on $n\ge n_0$ vertices.
Apply Theorem~\ref{kellythm} to $G$ in order to obtain a collection $\C$ of at least 
$(1/2-\eta_1)n$ edge-disjoint Hamilton cycles. Let $F$ be the undirected graph consisting of
all those edges of $G$ which are not covered by the Hamilton cycles in~$\C$. Note that $F$ is $k$-regular
for some $k\le 2\eta_1 n$. By Vizing's theorem the edges of $F$ can be coloured with at most
$\Delta(F)+1\le 3\eta_1 n$ colours and thus $F$ can be decomposed into at most $3\eta_1 n$ matchings.
Split each of these matchings into at most $1/\eta_1^{1/2}$ edge-disjoint matchings, each
containing at most $\eta_1^{1/2} n$ edges. So altogether this yields a collection $\M$ of
at most $3 \eta_1^{1/2} n$ matchings covering all edges of~$F$. It is enough to show that for
each $M\in \M$ there exists a Hamilton cycle of $G$ which contains all the edges in~$M$.

So consider any $M\in \M$. As observed in~\cite{KOTchvatal} (see the proof of Corollary~16 there),
any regular tournament is a robust $(\nu,\tau)$-outexpander. Let $D$ be the digraph obtained
from $G$ by `contracting' all the edges in~$M$, i.e.~by successively replacing each edge $xy\in M$ with a vertex
$v_{xy}$ whose inneighbourhood is the inneighbourhood of $x$ and whose outneighbourhood is
the outneighbourhood of $y$. Using that $M$ consists of at most $\eta_1^{1/2} n$ edges and that
$\eta_1\ll \nu,\tau$, it is not hard to check that $D$ is still a robust $(\nu/2,2\tau)$-outexpander
and $\delta^0(D)\ge (1/2-2\eta_1^{1/2}) n$. So Theorem~\ref{expanderthm} implies that
$D$ contains a Hamilton cycle, which corresponds to a Hamilton cycle in~$G$ containing all edges of~$M$,
as required. \endproof

Note that we cannot simply apply Theorem~\ref{main} instead of Theorem~\ref{expanderthm} at the end of the proof,
because $D$ may not been an oriented graph. However, 
instead of using Theorem~\ref{expanderthm}, one can also 
use the following result of Thomassen~\cite{thom2}: for every set $E$ of $n/24$
independent edges in a regular tournament on $n$ vertices, there is a Hamilton cycle
which contains all edges in $E$.

Theorem~\ref{RegularII} can be proved in a similar way, using Ghouila-Houri's theorem
instead of Theorem~\ref{expanderthm}.

\medskip
\noindent{\bf Proof of Theorem~\ref{RegularII}.} Choose a new constant $\eta$ such that
$1/n_0\ll\eta\ll \xi$ and apply Theorem~\ref{Regular} to find a collection of at least
$(d-\eta n)/2$ edge-disjoint Hamilton cycles. Let $F$ denote the subgraph of $G$ consisting 
of all edges not lying in these Hamilton cycles. Then $F$ is $k$-regular for some $k\le \eta n$.
Choose a collection $\M$ of matchings covering all edges of $F$ as in the the proof
of Theorem~\ref{cover}. So each matching consists
of at most $\eta^{1/2} n$ edges. As before, for each $M\in\M$ it suffices to find a Hamilton cycle of $G$
containing all edges of $M$. Let $D'$ be the digraph obtained from $G$ by orienting each
edge in $M$ and replacing each edge in $E(G)\setminus M$ with two edges, one in each direction.
Let $D$ be the digraph obtained from $D'$ by `contracting' the edges in $M$ as in the the proof
of Theorem~\ref{cover}. Then $D'$ has minimum semidegree at least $n/2$ and thus contains
a Hamilton cycle by Ghouila-Houri's theorem (Theorem~\ref{hamGH}). This Hamilton cycle
corresponds to a Hamilton cycle in~$G$ containing all edges of~$M$,
as required. \endproof

\section*{Acknowledgements} 
We are grateful to Demetres Christofides, Peter Keevash  and Andrew Treglown for their comments on a 
draft of this paper.

{\footnotesize
\bigskip\obeylines\parindent=0pt
Daniela K\"uhn \& Deryk Osthus
School of Mathematics
Birmingham University
Edgbaston
Birmingham B15 2TT
UK
{\it E-mail addresses}: {\tt \{kuehn,osthus\}@maths.bham.ac.uk}
}

\end{document}